# Deep Ritz Physics-Informed Neural Network Method for Solving the Variational Inequality Problems

Qijia Zhou, Yiyang Wang, Shengyuan Deng, Chenliang Li *

**Abstract**—Variational inequalities are widely applied in mechanical engineering, fluid penetration, transportation, and other fields. In this paper, a Deep Ritz method based on Physics-Informed Neural Networks (PINNs) is proposed to enhance the accuracy and efficiency of solving elliptic variational inequalities. The Ritz variational method is firstly utilized to transform the variational inequality problem into an optimization problem. Then Bayesian optimization is employed to tune the weights of the loss function, and a residual-based adaptive dataset update strategy is introduced to improve the convergence and accuracy of the model. Numerical experiments show that the proposed method can effectively approximate the analytical solution.

## I. INTRODUCTION

Elliptic variational inequalities and complementarity problems are widely applied in many fields, for example, mechanical engineering, fluid penetration, transportation and others. In recent years, many significant numerical algorithms have been presented to solve the variational inequality and the complementarity problems, but often involving expensive computation cost and the complexity.

The artificial neural network has the abilities of massively parallel processing and good stability. Therefore, many neural network models have been presented for solving the variational inequalities. Some neural-networks are presented for the variational inequalities in [1-3], a novel neural network is discussed for variational inequalities with linear and nonlinear constraints in [4], some deep learning-based numerical algorithms are given for solving variational problems in [5-7].

Physics-Informed Neural Networks (PINNs) do not merely recognize distribution trends in learning datasets like conventional neural models; they also encapsulate fundamental principles of physics. By embedding physical constraints during the learning phase, PINNs generate more generalized models with reduced data dependency[8, 9]. Currently, PINN models have been extensively applied across diverse domains. For example, PINN-FORM integrates PINNs with FORM to enhance the reliability assessment of partial differential formulations, thereby improving predictive accuracy[10,11]. Additionally, a self-adaptive PINN-based framework has been proposed to address nonlinear partial differential equations efficiently[12]. Furthermore, a novel phase-field smoothing approach[13], leveraging physics-aware neural architectures, has been introduced to solve differential equations featuring discontinuous coefficients. Experimental results confirm that these methodologies exhibit high accuracy and efficacy in tackling elliptic equations. Deep learning methods based on PINNs have been widely concerned, such as deep domain decomposition method[14], PINNs for fluid mechanics[15] and wave equations[16], fast PINN network[17], Self-adaptive PINNs[18] and error estimate of PINNs to approximating PDEs[19], and so on.

However, the application of PINNs in elliptic variational inequalities and complementarity problems has not been widely explored. Therefore, we try to construct a deep learning method based on PINNs in this field. Firstly, we theoretically convert the variational inequalities into an optimization problem. Subsequently, for enhancing the precision of the solution, we introduce Bayesian optimization and adaptive residual methods. Finally, through simulation and comparison with other methods, we aim to assess the accuracy of the proposed approach.

## II. DEEP RITZ-PINNS FOR VARIATIONAL INEQUALITIES

### A. Ritz Method

In this paper, we consider the elliptic variational inequality problem (EVIP): find $u \in V = \{v \in H_0^1(\Omega), v \geq 0 \ a.e. \ \Omega\}$, such that

$$v \in V, \quad a(u, v-u) \geq f(v-u), \quad (1)$$

where

$$a(u,v) = \int_\Omega (\nabla u \cdot \nabla v + \alpha uv) d\Omega, \quad f \in L^2(\Omega), \quad (2)$$

and $\Omega$ is a polygon.

Based on the Ritz method, it is evident that problem (1) can be reformulated as the following optimization model, finding $u \in V$ such that

$$J(u) = min J(v), \quad v \in V,$$

Manuscript received March 31, 2025; revised October 11, 2025.
This work is supported by the National Undergraduate Training Program for Innovation (202310595022).

Qijia Zhou is an undergraduate student of School of Mathematics and Computing Science, Guilin University of Electronics Technology, Center for Applied Mathematics of Guangxi (GUET), Guilin 541004, Guangxi, China (email: denglong836@163.com).

Yiyang Wang is an undergraduate student of School of Mathematics and Computing Science, Guilin University of Electronics Technology, Center for Applied Mathematics of Guangxi (GUET), Guilin 541004, Guangxi, China (email: gongyou6@163.com).

Shengyuan Deng is an undergraduate student of School of Mathematics and Computing Science, Guilin University of Electronics Technology, Center for Applied Mathematics of Guangxi (GUET), Guilin 541004, Guangxi, China (email: sarengchenguv41@163.com).

Chenliang Li is a Professor of School of Mathematics and Computing Science, Guilin University of Electronics Technology, Center for Applied Mathematics of Guangxi (GUET), Guilin 541004, Guangxi, China (Corresponding author's email: laiyingg6793@163.com).

where
$$J(v) = \frac{1}{2}a(v,v) - f(v).$$

Another way, problem (e1) is also equivalent to the following complementarity problem: find $u \in V$ such that
$$\begin{cases} -\Delta u + \alpha u - f \geq 0 \\ u(-\Delta u + \alpha u - f) = 0 \\ u \geq 0 \end{cases}$$

*B. Physics-Informed Neural Networks (PINNs)*

PINNs is an approach for solving Partial Differential Equations (PDEs) by embedding physical laws into the neural network optimization process. This is achieved by incorporating the governing PDEs and boundary conditions into the loss function, guiding the network to respect physical constraints. Therefore, it is particularly effective when dealing with complex physical problems where obtaining large amounts of experimental data is difficult. We construct a class of PINNs to solve problem (2).

A classical PINN network has $L + 1$ layers and the Layer 0 is the input layer. The Layer L is the output layer, and Layer 1 to Layer L are the hidden layers. Each hidden layer contains an activation function, such as Sigmoid, Tanh, or ReLU, while the output layer has no activation function. The activation function introduces nonlinearity into the network, enabling it to approximate complex solutions to nonlinear PDEs. A detailed introduction of each layer is presented below.

**Input layer**: The input data is represented as a vector $x \in R^n$, n denotes the dimension of the input data. This vectorized form allows the neural network to process inputs in a unified mathematical structure.

**Hidden layer**: These layers are calculated through a linear transformation and activation function. Let the $l^{th}$ layer ($l = 1, 2, 3, \ldots, L$) have $M_l$ neurons, and the parameter space is defined as $\theta := \{W_l, b_l\}$ represents the collection of all parameters, Where $W_l$ represents the weights and $b_l$ represents the biases. Define function $f_l: \mathbb{R}^{d_l} \to \mathbb{R}^{d_{l+1}}, 0 \leq l < L$
$$f_l = W_l * \sigma(x) + b_l$$

Where $\sigma(x)$ is the activation function. For example, we can choose ReLU function as the activation function.

The operations of all L hidden layers can be represented as a composite function, resulting in the composite function for the entire hidden layer part
$$f_{hidden}(x) = f_L \circ f_{L-1} \circ \cdots \circ f_0$$

**Output layer:** The activation function in the output layer is an identity function,
$$\alpha(z) = z$$

**Loss function:** The objective function $M(\theta)$ is defined as follows
$$\theta^* = \text{argmin}_\theta M(\theta)$$
$$:= w_1 M_\Omega(\theta) + w_2 M^+{}_\Omega(\theta) + w_3 M_{\partial\Omega}(\theta),$$
where
$$M_\Omega(\theta) := \frac{1}{N_f}\sum_{i=1}^{N_f}\left(\frac{1}{2}|\nabla_x u(\mathbf{x}_f^i;\theta)|^2 - f(\mathbf{x}_f^i)u(\mathbf{x}_f^i;\theta)\right), \quad (3)$$
$$M^+(\theta) := \frac{1}{N_f}\sum_{i=1}^{N_f} \max\{-u(\mathbf{x}_f^i;\theta), 0\}, \quad (4)$$
$$M_{\partial\Omega}(\theta) := \frac{1}{N_g}\sum_{i=1}^{N_g}|u(\mathbf{x}_f^i;\theta) - g(\mathbf{x}_g^i)|^2 \quad (5)$$

$\{\mathbf{x}_f^i\}_{i=1}^{N_f}$ and $\{\mathbf{x}_g^i\}_{i=1}^{N_g}$ are the collocation points. The domain term $M_\Omega$ and boundary term $\mathcal{M}_{\partial\Omega}$ enforce the condition that the desired optimized neural network $h(;\theta^*)$ satisfes $a(u,u) = f(u)$ and $\mathcal{B}(u) = g$, $w_1, w_2$ and $w_3$ are the weights of loss function.

*C. Adam and Bayesian Optimization*

After constructing PINN model, we need to update and optimize the parameters to minimize the loss function and ensure that the model can accurately learn the solution. In this process, we choose to use the Adam optimizer (Adaptive Moment Estimation), which is a first-order gradient-based optimization method. The Adam optimizer, integrating momentum and adaptive learning rate techniques, makes it well-suited for handling sparse gradients and non-stationary objectives. It is particularly well-suited for large-scale data and high-dimensional parameter optimization. Therefore, choosing Adam allows the PINN model to accelerate training and improve accuracy when dealing with complex physical problems.

---

**Algorithm 1**: Adam

**Require**: $\alpha$: Stepsize
**Require**: $\beta_1, \beta_2 \in [0,1)$, which are used to control momentum term and RMSprop-like term
Require: $M(\theta)$, which is the stochastic objective function
**Require**: $\theta_0$: Initial parameter vector, let $m_0$, $v_0$ and $t \leftarrow 0$
While(epoch)
   while $\theta_t$ not converged do
      1. $t \leftarrow t + 1$
      2. $g_t \leftarrow \nabla_\theta M_t(\theta_{t-1})$
      3. $m_t \leftarrow \beta_1 \cdot m_{t-1} + (1 - \beta_1) \cdot g_t$
      4. $v_t \leftarrow \beta_2 \cdot v_{t-1} + (1 - \beta_2) \cdot g_t^2$
      5. $\widehat{m}_t \leftarrow m_t/(1 - \beta_1^t)$
      6. $\hat{v}_t \leftarrow v_t/(1 - \beta_2^t)$
      7. $\theta_t \leftarrow \theta_{t-1} - \alpha \cdot \widehat{m}_t/(\sqrt{\hat{v}_t} + \epsilon)$
   end while
   return $\theta_t$
end

Parameters description:
1. $t$: Number of steps (steps)
2. epoch: Number of training iterations.
3. $\alpha$: learning rate, Used to control the step size.
4. $\theta$: Parameters to be solved (updated).
5. $M(\theta)$: The stochastic objective function with parameters $\theta$, generally referring to the loss function.
6. $g_t$: The gradient of $f(\theta)$ with respect to $\theta$ is obtained by taking the derivative.
7. $\beta_1$: First-order moment decay coefficient.
8. $\beta_2$: Second-order moment decay coefficient.
9. $m_t$: The first-order moment of the gradient $g_t$, representing its expected value.
10. $v_t$: The first-order moment of the gradient $g_t$, also representing its expected value.
11. $\widehat{m}_t$: The bias correction of $m_t$.
12. $\widehat{v}_t$: The bias correction of $v_t$.

---

For the selection of weights, we optimize them by using the Bayesian optimization method to achieve a balance between the error terms of different loss functions. Compared to manually tuning parameters, the Bayesian optimization method can infer the optimal values of the weights based on observational data, balancing efficient exploration and exploitation, thus maximizing the optimization efficiency of the loss function. Typically, when there are multiple weights, people tend to try different combinations one by one, and the choice of each weight can have a different impact on the result. However, Bayesian optimization can effectively assign appropriate weights to improve the overall optimization process. The algorithm is presented as follows.

---

**Algorithm 2:** Bayesian Optimization for Weight

**Initialization:** Loss function $loss(W) = w_1 * loss_1 + w_2 * loss_2 + w_3 * loss_3$,
where
$W = \{w_1, w_2, w_3\}, \ D = \{(w_i, loss(w_i))\}$.
**Require**:
  1. Optimize the function $loss(W)$ with the set of weights $w_{new}$
  2. Evaluate the loss function $loss(w_{new})$
  3. Update the dataset
$$D = D \cup \{w_{new}, loss(w_{new})\}$$
**Return**: The normalized optimized weights $w_1, w_2, w_3$.

---

*D. Residual-based Dataset Update*

Different from the traditional training methods based on a fixed dataset, we introduce a residual-adaptive dataset update strategy. In the traditional training process, the model adjusts the weights to focus on key directions. In the residual-adaptive training process, the model dynamically focuses on areas with larger residuals (prediction errors), allowing the model to concentrate on difficult-to-fit regions. Additionally, during the training procedure, points in the dataset can be dynamically adjusted, gradually guiding the model to learn more complex regions. With this strategy, the model can more efficiently capture the important features in the problem, thereby improving the overall training effectiveness and solution accuracy. The specific steps are as follows.

**Step 1: Initial sampling**

In the domain $\Omega$, the initial point set $D_0 = \{(x_i, y_i) | i = 1, 2, \dots, N\}$ is generated using Latin Hypercube Sampling (LHS).

**Step 2: Calculate the residual**

In each training iteration, the current network prediction $u(x, y)$ is substituted into the loss function to calculate the residual for all sampling points
$$R(x, y) = loss(x, y)$$

**Step 3: Calculate the sampling weights**.

Normalize the residual $R(x, y)$, and obtain the sampling probability.
$$P(x_i, y_i) = \frac{R(x_i, y_i) + \varepsilon}{\sum R(x', y') + \varepsilon}$$

The normalization here ensures that the sampling probability of all points sums to 1, and $\varepsilon = 10^{-6}$ is used to avoid numerical issues.

**Step 4: Dynamically generate new sampling points**

M parent points are selected through sampling with replacement. For each parent point $(x_p, y_p)$, generate n child points within its domain. The radius $r = r_0 \times e^{-t}$, centered at the parent point, is adaptively adjusted, where $t$ is the iteration count and $r_0$ is the initially set radius.

**Step 5: Update the sampling point set**

The newly generated sampling points replace the previous point set, where
$$D_{t+1} = D_t \cup D_{new} - D_{point}$$

To prevent duplicate points in the dataset, the previously selected centroid points $D_{point}$ are removed.

**Step 6: Re-training**

Continue training using the updated sampling point set $D_{t+1}$.

*E. Algorithmic Implementation of the Deep Ritz-PINNs*

In summary, based on the Adam, Bayesian optimization for weights, and residual-based dataset update method, we propose a new method, called Deep Ritz PINNs method, to efficiently solve the variational inequality problems. The algorithm steps are as follows.

---

**Algorithm 3:** Deep Ritz-PINNs (DRPINNS)

While ($epoch$):
  1. Train PINN using the Adam optimizer and current data points
   - Compute the total loss,
      $loss = w_1 * loss_1 + w_2 * loss_2 + w_3 * loss_3$
   - Update network parameters by minimizing the loss using Adam steps
  2. Optimize the loss weights $w_1, w_2, w_3$ by using Bayesian Optimization
   - Fit a surrogate model (e.g., Gaussian Process) to predict loss values from observed weights
   - Propose new weights that minimize the predicted loss, by maximizing an acquisition function (e.g., Expected Improvement).
   - Update the Bayesian optimization dataset with the new sample point.
  3. Update the training dataset (if adaptive sampling is applied or constraint is triggered)
  4. Check stopping criteria (maximum number of epochs or loss threshold).
end while

Return: the approximate solution $u(x, y)$,.

---

The algorithm is shown in Figure 1.

III. NUMERICAL EXAMPLES

In this section, the network architecture is simply set with the same number of units in each layer. The "number of layers" refers to the hidden layers in the network, while the "number of units" refers to the number of units in each layer. We choose Tanh as the activation function and employ the Adam optimization algorithm for updating parameters via stochastic gradient descent. In this paper, the number of iterations is set to 100k or until the error between iterations becomes less than $10^{-4}$. The learning rate is adjusted based on different conditions, with the initial learning rate and decay rate being modifiable. Specifically, the initial learning rate ranges from 0.0001 to 0.01. The initial training data in this paper is randomly generated according to a normal distribution, and the dataset is updated every 10k iterations or when the iteration error is less than $10^{-3}$. The initial weights are set to $10^4$ to address the issue of vanishing gradients in the model.

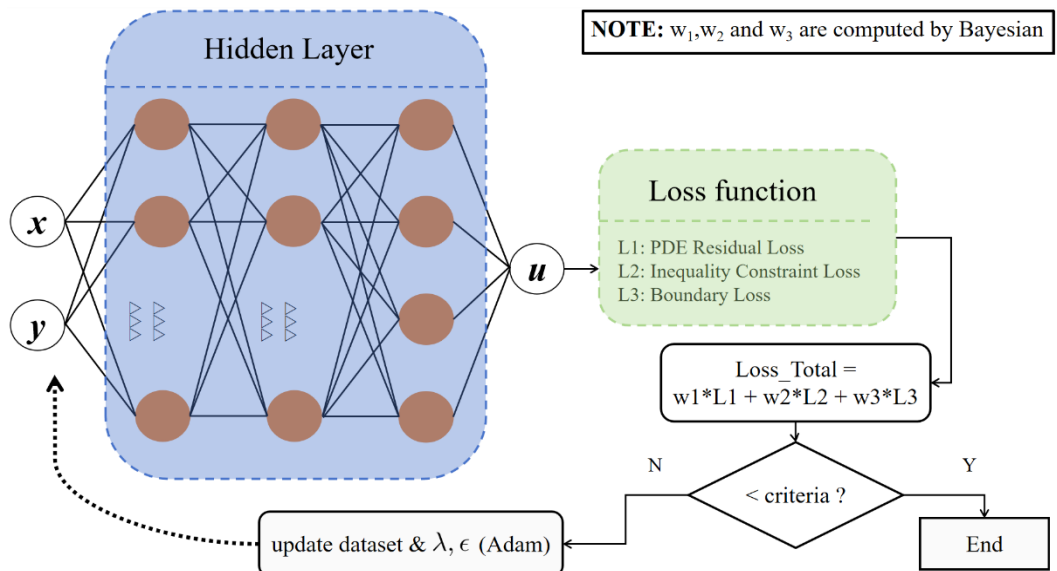

Figure 1 DRPINNs network

Additionally, the computational setup used for this study includes a CPU with 192 cores running at a frequency of 1130.5 MHz, along with an RTX 4090D (24GB) graphics processor, which facilitates the high-performance training and optimization processes.

*A. 1D Obstacle Problem*

**Example 1** Given the following variational inequality, finding $u \in K$, such that

$$\begin{cases} \int_0^1 u'(x)(v'(x) - u'(x))dx \geq \int_0^1 f(x)(v(x) - u(x))dx, \forall v \in K, \\ u(x) \geq \psi(x), x \in (0,1), \\ u(0) = u(1) = 0, \end{cases}$$

where $f(x) = 2x$ and $\psi(x) = \frac{1}{4}x(1-x)$.

First, expand the integral on the left-hand side of the inequality.

$$\int_0^1 u'(x)v'(x)dx - \int_0^1 (u'(x))^2 dx$$
$$\geq \int_0^1 f(x)(v(x) - u(x))dx.$$

After rearrangement, we obtain

$$\int_0^1 u'(x)v'(x)dx - \int_0^1 f(x)v(x)dx$$
$$\geq \int_0^1 (u'(x))^2 dx - \int_0^1 f(x)u(x)dx.$$

Since $u$ satisfies the boundary conditions, the problem reduces to finding $u$ that minimizes the following functional

$$J(v) = \frac{1}{2}\int_0^1 (v'(x))^2 dx - \int_0^1 f(x)v(x)\,dx.$$

Therefore, the equivalent minimization problem can be formulated as finding $u \in k$, such that

$$J(u) = \min_{v \in K} J(v) \quad (6)$$

$$\text{s.t.} \begin{cases} u(x) \geq \psi(x), x \in (0,1), \\ u(0) = u(1) = 0. \end{cases}$$

After transforming the above equation into the loss function, the approximate solution of the model is computed through DRPINNs. First, plot the graphs of the exact solution and the analytical solution.

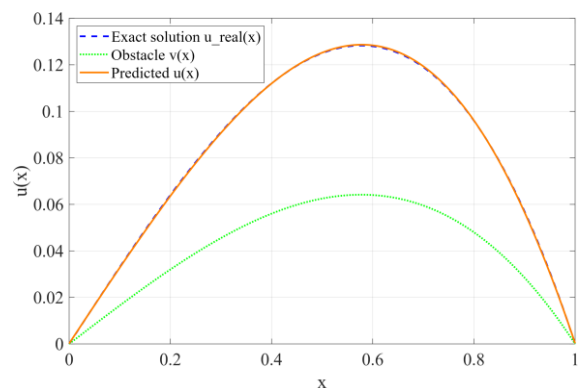

Figure 2 Exact solution vs predicted solution of Example 1.

As shown in Figure 2, the predicted solution $u(x)$ maintains a high degree of consistency with the exact solution $u_{real}(x)$ across the entire domain, while strictly satisfying the obstacle constraint $u(x) \geq \psi(x)$. Moreover, the predicted solution successfully exhibits the characteristics of "adhering to" and "detaching from" the obstacle function near the obstacle region, demonstrating the effectiveness and accuracy of the proposed method in handling one-dimensional variational inequality problems.

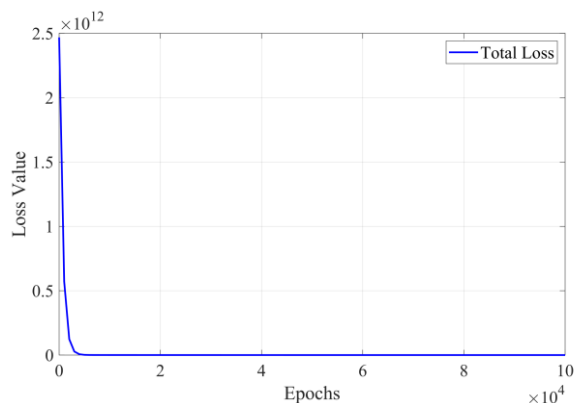

Figure 3 Loss function in Example 1

During the training process, the total loss decreases rapidly and eventually converges to a stable value, as shown in Figure 3. This indicates that the model effectively fits the problem and that the training process exhibits good convergence and stability.

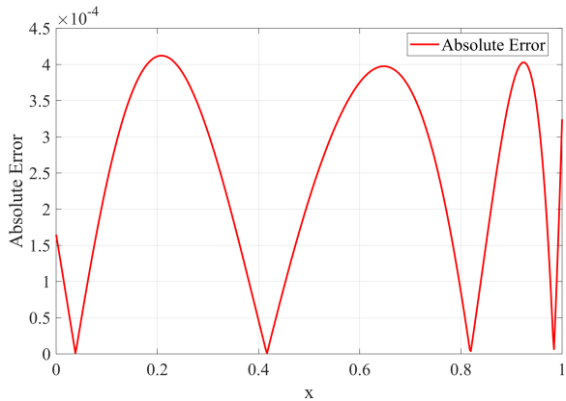

Figure 4 Absolute error in Example 1

Figure 4 shows the distribution of the absolute error between the predicted solution and the exact solution. It can be observed that the overall magnitude of the error remains small, with fluctuations in several local regions, demonstrating good global fitting accuracy and stability of the model.

*B.2D Elliptic Variational Inequality*

**Example 2** Given the following variational inequality, finding u, such that
$$\begin{cases} -u_{xx} - u_{yy} \geq 2\pi^2 \sin\pi x \sin\pi y, (x,y) \in \Omega = [-1,1] \times [-1,1] \\ u \geq 0, \\ u(-u_{xx} - u_{yy} - 2\pi^2 \sin\pi x \sin\pi y) = 0, \\ u = 0, on\ \partial\Omega. \end{cases}$$

Firstly, by multiplying both sides of the inequality by a test function $v(x,y)$ and integrating over the domain $\Omega$, we obtain
$$-\int_\Omega (u_{xx} + u_{yy}) v dx dy \geq \int_\Omega 2\pi^2 \sin(\pi x) \sin(\pi y)\ v dx dy,$$

According to Green's formula, we have
$$-\int_\Omega (u_{xx} + u_{yy}) v dx dy = \int_\Omega \nabla u \cdot \nabla v dx dy - \int_{\partial\Omega} \frac{\partial u}{\partial n} v ds,$$

Because $u = 0$ on $\partial\Omega$, we have
$$-\int_\Omega (u_{xx} + u_{yy}) v dx dy = \int_\Omega \nabla u \cdot \nabla v dx dy$$

Substituting Green's formula into the original inequality
$$\int_\Omega \nabla u \cdot \nabla v dx dy \geq \int_\Omega 2\pi^2 \sin(\pi x) \sin(\pi y)\ v dx dy,$$

Therefore, the equivalent minimization problems is as follows, finding $u \in V = \{v \in H_0^1(\Omega), v \geq 0\ a.e.\ \Omega\}$, such that
$$J(u) = \min J(v)$$
$$= \frac{1}{2} \int_0^1 \int_0^1 [(v_x)^2 + (v_y)^2] dx dy - \int_0^1 \int_0^1 2\pi^2 \sin(\pi x) \sin(\pi y) \cdot v dx dy.$$

The final expression is given by
$$J(u) = \min J(v) \qquad (7)$$

$$\begin{cases} u \geq 0, \\ u(-u_{xx} - u_{yy} - 2\pi^2 \sin\pi x \sin\pi y) = 0, \\ u = 0, on\ \partial\Omega. \end{cases}$$

Following the same procedure as before, we first plot the predicted solution and the analytical solution.

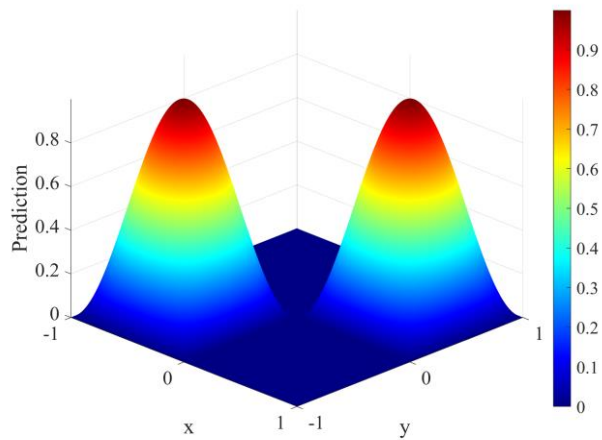

Figure 5 Predicted solution of Example 2

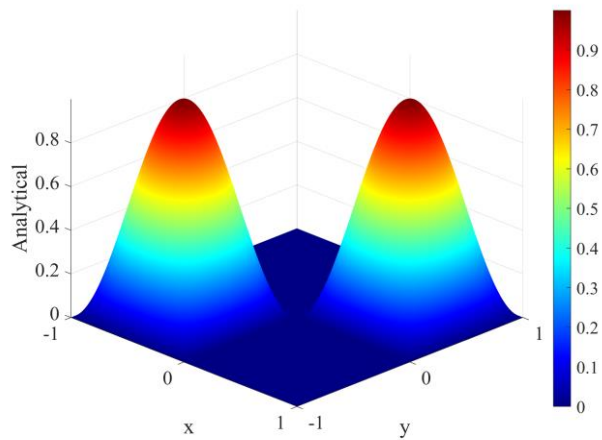

Figure 6 Analytical solution of Example 2

Figure 5 shows the numerical solution computed by DRPINNs, while Figure 6 shows the analytical solution of the problem. As can be seen, the model's output is very close to the analytical solution of the original problem, indicates that the training of the model has been successful. Next, we plot both the 3D error and the 2D error graphs to examine the errors in detail.

As can be seen from Figure 7, the error distribution is closely related to the local characteristics of the solution function in different regions. Specifically, the solution changes drastically in the 1st and 3rd quadrant regions, and the model fitting is difficult, so the corresponding error is high. In the second and fourth quadrants, the function value tends to be zero, the whole is relatively flat, and the error is reduced. Overall, the maximum error is about 0.001, which indicates that the model has good prediction accuracy. However, there are still some fluctuations in the error graph, which indicates that there is still room for improvement in the expressive ability of the model when dealing with the high gradient region. This fluctuation may be due to the difference in regularity of the solution in different regions and the non-uniform error distribution caused by the complex boundary structure.

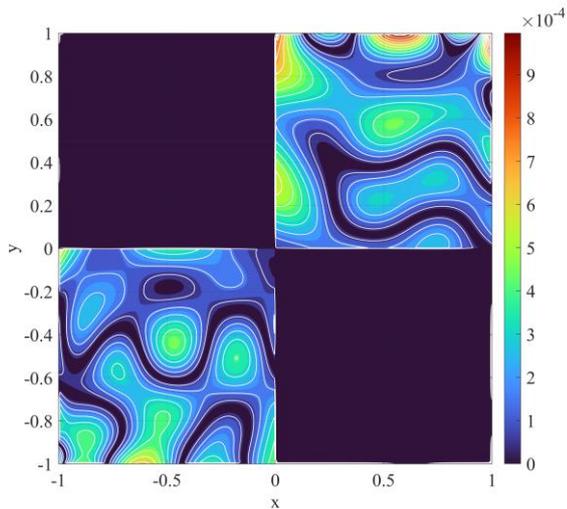

Figure 7 2D error in Example 2

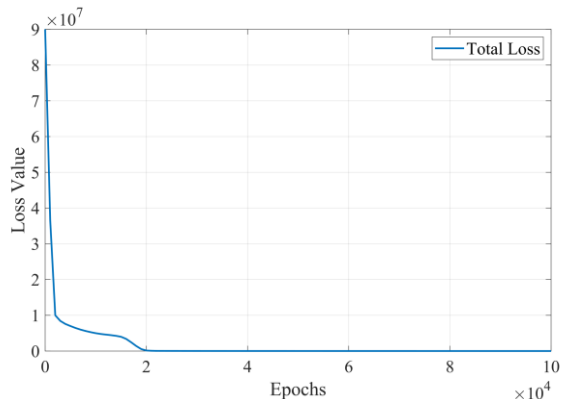

Figure 8 Loss function in Example 2

Figure 8 shows the loss function. Through analysis of the total loss function, it can be observed that in the early stages of training, the total loss rapidly decreases, which indicates that the model learns a significant amount of useful information from the data during the initial training phase. Subsequently, the descent trend of the loss function gradually becomes gentler, which indicates that the model has converged, and meaning that it has learned the features of the data and reached a stable state. To observe the changes in the loss function more clearly, we record the loss function every 5k iterations and set different epochs intervals for separate viewing.

*C. 2D Elliptic Variational Inequality with Spatially Dependent Source Term*

**Example 3** Given the following variational inequality, finding u, such that

$$\begin{cases} -\Delta u \geq 4 - 2(x^2 + y^2), & (x,y) \in \Omega = [-1,1] \times [-1,1] \\ u \geq 0, \\ u(-\Delta u - 4 + 2(x^2 + y^2)) = 0, \\ u = 0, & on\ \partial\Omega \end{cases}$$

Firstly, by multiplying both sides of the inequality by a test function $v(x,y)$ and integrating over the domain $\Omega$, we obtain

$$-\int_\Omega (\Delta u) v\ dxdy \geq \int_\Omega (4 - 2(x^2 + y)) v\ dxdy.$$

According to Green's formula, we have

$$-\int_\Omega (\Delta u) v\ dxdy = \int_\Omega \nabla u \cdot \nabla v\ dxdy - \int_{\partial\Omega} \frac{\partial u}{\partial n} v\ ds.$$

Because $u = 0\ on\ \partial\Omega$ and $u = 0\ on\ \partial\Omega$, we have

$$-\int_\Omega (\Delta u) v\ dxdy = \int_\Omega \nabla u \cdot \nabla v\ dxdy.$$

Substituting Green's formula into the original inequality, we get

$$\int_\Omega \nabla u \cdot \nabla (v - u)\ dxdy$$
$$\geq \int_\Omega (4 - 2(x^2 + y^2))(v - u)\ dxdy.$$

Therefore, the equivalent minimization problem is as follows, finding $u \in V = \{v \in H_0^1(\Omega), v \geq 0\ a.e.\ \Omega\}$, such that

$$J(u) = \min J(v) = J(v)$$
$$= \frac{1}{2} \int_\Omega |\nabla v|^2\ dxdy - \int_\Omega (4 - 2(x^2 + y^2))\ v\ dxdy.$$

The final expression is given by

$$J(u) = \min J(v) \qquad (8)$$
$$s.t. \begin{cases} u \geq 0, \\ u(-\Delta u - 4 + 2(x^2 + y^2)) = 0, \\ u = 0, \quad on\ \partial\Omega. \end{cases}$$

Based on (8), we plot the predicted solution and the analytical solution.

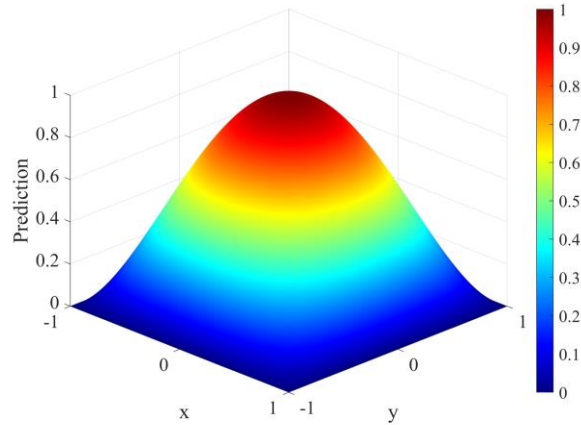

Figure 9 Predicted solution of Example 3

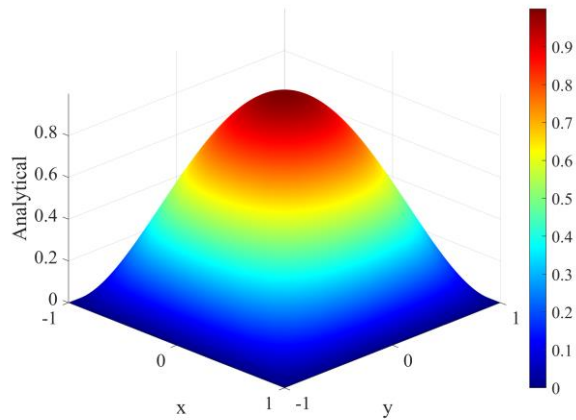

Figure 10 Analytical solution of Example 3

The predicted solution closely matches the exact solution, demonstrating the accuracy of the method.

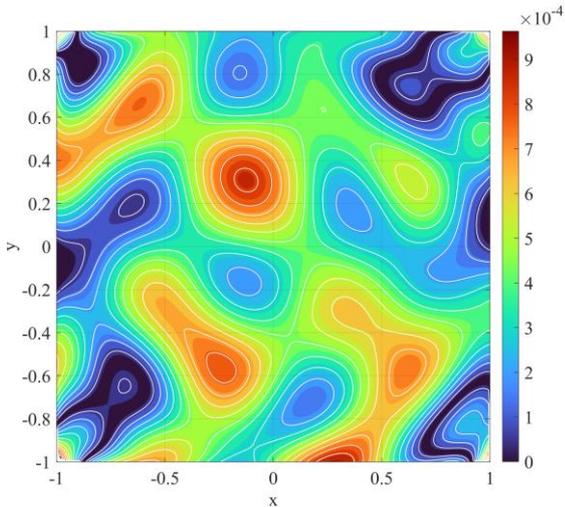

Figure 11 2D error in Example 3

From Figure 11, it can be observed that the error is uniformly distributed with small magnitudes, with only slight fluctuations in some local regions, which indicates the good global approximation accuracy and stability of the model.

**Example 4** Given the following variational inequality, finding u, such that
$$\begin{cases} -\Delta u \geq f(x,y,z), (x,y,z) \in \Omega = [0,1] \times [0,1] \times [0,1] \\ u \geq 0, \\ u(-\Delta u - f(x,y,z)) = 0, \\ u = r^2 - r_0^2, \quad on \ \partial\Omega \end{cases}$$

where
$$f(x,y,z) = \begin{cases} -4(2r^2 + 3(r^2 - r_0^2)), & r > r_0, \\ -8r_0^2(1 - r^2 + r_0^2), & r \leq r_0, \end{cases}$$
where $r = (x^2 + y^2 + z^2)^{1/2}$ and $r_0 = 0.7$.

The non-homogeneous Dirichlet boundary condition is taken in such a way that the solution $u$ is given by $u = (max(r^2 - r_0^2, 0))^2$.

Firstly, the variational form is as follows,
$$\int_\Omega \nabla u \cdot \nabla(v - u) dx \geq \int_\Omega f(x,y,z)(v - u)dx, \forall v \geq 0.$$

According to Green's formula, Combined with non-homogeneous Dirichlet boundary conditions $(u|_{\partial\Omega} = r^2 - r_0^2)$

Therefore, the equivalent optimization problem is
$$u = \arg \min_{u \in K} J(u),$$
where
$$J(u) = \frac{1}{2} \int_\Omega |\nabla u|^2 dx - \int_\Omega f(x,y,z) u dx,$$
$$K = \{u \in H^1(\Omega) \mid u \geq 0, u|_{\partial\Omega} = r^2 - r_0^2\}.$$
The final expression is given by

$$J(u) = \min J(v)$$
$$s.t. \begin{cases} u \geq 0, \\ u(-\Delta u - 4 + 2(x^2 + y^2)) = 0, \\ u = r^2 - r_0^2, \quad on \ \partial\Omega. \end{cases}$$

To provide a clearer illustration of the solution, we compare the analytical and predicted results. Since it is challenging to directly visualize the solution in three dimensions, we transform the Cartesian coordinates $(x, y, z)$ into the radial coordinate
$$r = \sqrt{x^2 + y^2 + z^2}$$
and plot the solution as a function of $r$. It can be observed that the solution remains zero for $r < r0 = 0.7$, while for $r \geq r0$, the predicted solution closely matches the analytical one, with the two curves nearly overlapping. The final result is shown in Figure 12.

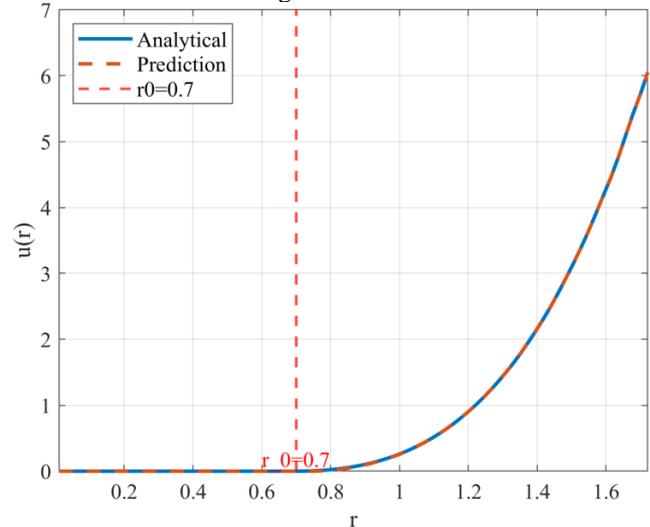

Figure 12 Exact solution and predicted solution of Example 4

The relative max error of the solution is shown in Figure 13.

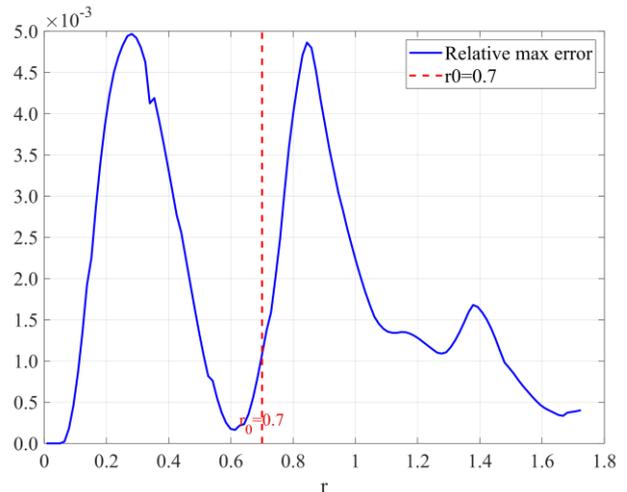

Figure 13 Relative max error of Example 4

*D. Ablation Experiments*

In this section, we present a systematic analysis of the DRPINNs method through ablation studies. The goal is to quantitatively evaluate the impact of key mechanism (Bayesian optimization and dataset update strategies) on model performance, and to verify their effectiveness in error control and training efficiency improvement. To this end, we compare the full DRPINNs method with its variants: one without the Bayesian optimization module and one without the dataset update strategy. The evaluation is conducted based on four metrics.

MSE stands for Mean Squared Error as followings,
$$MSE = \frac{1}{n} \sum_{i=1}^{n} (y_i - \hat{y}_i)^2$$

MAE stands for Mean Absolute Error,

$$MAE = \frac{1}{n}\sum_{i=1}^{n}|y_i - \hat{y}_i|$$

REL2 stands for Relative $L2$ Norm Error which is defined as,

$$REL2 = \frac{\|y - \hat{y}\|_2}{\|y\|_2}$$

and MAXERROR is defined by the infinity norm as

$$MAXERROR = \max_i |y_i|$$

The relative maximum error is defined as

$$\text{Relative-MAXERROR} = \frac{\max_i |y^i - \hat{y}^i|}{\max|y|}.$$

Here, n represents the number of samples, $y^i$ denotes the actual observations, and $\hat{y}^i$ represents the corresponding predicted values.

The experimental results are presented in the following Table 1.

Table 1 COMPARISON OF EVALUATION METRICS

| | Method | MSE | MAE | REL2 | (relative) MAXERROR |
|---|---|---|---|---|---|
| Exp1 | DRPINNS | 3.06e-8 | 1.02e-4 | 2.49e-7 | 9.04e-4 |
| | DRPINNs without dataset update | 7.23e-7 | 4.97e-4 | 5.90e-6 | 2.90e-3 |
| | DRPINNs without Bayesian optimization | 9.25e-2 | 2.83e-1 | 7.55e-1 | 5.97e-1 |
| Exp2 | DRPINNS | 7.91e-8 | 2.51e-4 | 3.06e-3 | 4.12e-4 |
| | DRPINNs without dataset update | 2.32e-3 | 4.23e-2 | 5.24e-1 | 7.25e-2 |
| | DRPINNs without Bayesian optimization | 4.44e-6 | 1.76e-3 | 2.29e-2 | 3.64e-3 |
| Exp3 | DRPINNS | 1.87e-7 | 3.84e-4 | 6.72e-7 | 9.82e-4 |
| | DRPINNs without dataset update | 7.10e-7 | 6.29e-4 | 2.55e-6 | 2.75e-3 |
| | DRPINNs without Bayesian optimization | 1.42e-1 | 3.73e-1 | 5.09e-1 | 4.75e-1 |
| Exp4 | DRPINNS | 4.66e-5 | 4.01e-3 | 7.46e-3 | 5.62e-3 |
| | DRPINNs without dataset update | 7.41e-5 | 6.01e-3 | 9.41e-3 | 7.68e-3 |
| | DRPINNs without Bayesian optimization | 1.29 | 1.00 | 1.24 | 7.49e-1 |

From Table 1, we know that, with dataset updating and Bayesian optimization, DRPINNs achieves the lowest errors, The results indicate that, dynamically updating the training dataset contributes to improved accuracy and generalization, and Bayesian optimization plays a crucial role in enhancing the training efficiency and achieving better convergence.

*E. Comparative Experiment*

To validate the effectiveness of the proposed method (DRPINNs), we selected two representative deep learning approaches for comparison: Deep Equilibrium Models (DEQ) [20] and Augmented Lagrangian Deep Learning (ALDL) [21]. The DEQ method replaces the traditional stacked network architecture by solving for the equilibrium point of a neural network, thereby achieving constant memory cost and an effectively infinite-depth structure. The ALDL method, on the other hand, is based on the augmented Lagrangian framework, employing two neural networks to approximate the primal variables and the Lagrange multipliers separately, enabling effective handling of variational problems with essential boundary conditions. To ensure a fair comparison, all methods use the same hyperparameter settings: a neural network with three hidden layers, a fixed number of nodes per layer, 2000 training points, and a learning rate of 1e-4, without additional hyperparameter tuning for individual methods. The following sections present the error evolution curves of the three methods on the test problems and compare their convergence speeds and final error levels.

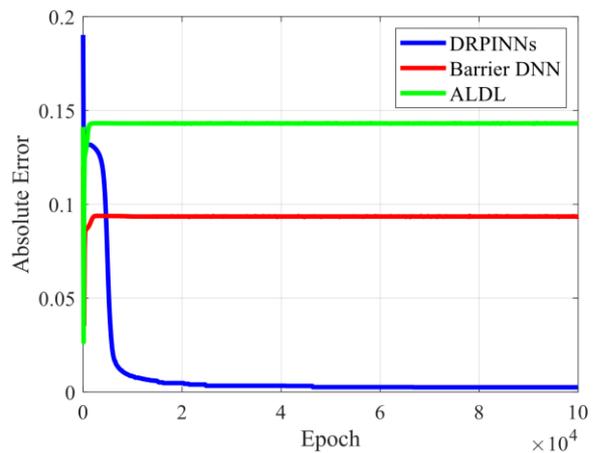

Figure 14 Error comparison for DRPINNs, Barrier DNN, and ALDL on Example 1

For the one-dimensional obstacle variational inequality problem (Example 1), we compare the training processes and error evolutions of three different methods. As shown in Figure 14, the DRPINNs method exhibits a rapid decrease in error from the early stages of training and eventually stabilizes at an extremely low error level, demonstrating excellent convergence and approximation capabilities for constrained problems. In contrast, the Barrier DNN method stagnates after an initial phase of convergence and results in a relatively large final error. The ALDL method, on the other hand, shows significant oscillations throughout the training process, with poor error control. Overall, DRPINNs achieves the best numerical performance for this problem.

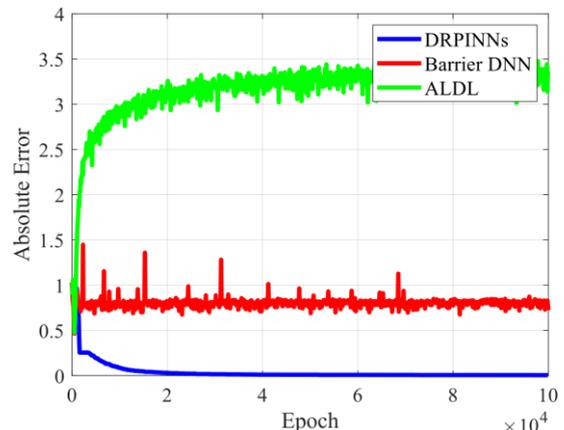

Figure 15 Error comparison for DRPINNs, Barrier DNN, and ALDL on Example 2

In Example 2, the DRPINNs method once again demonstrates outstanding performance, as shown in Figure 15. Within a short period, its absolute error rapidly decreases to an extremely low level and remains stable thereafter. In contrast, the Barrier DNN method converges more slowly and exhibits considerable fluctuations during the training process, resulting in a final error much higher than that of DRPINNs. Meanwhile, the ALDL method shows a trend of continuous deterioration throughout training, with the error steadily increasing and failing to converge effectively. Overall, DRPINNs significantly outperforms the two comparative methods in terms of convergence speed, stability, and final error control.

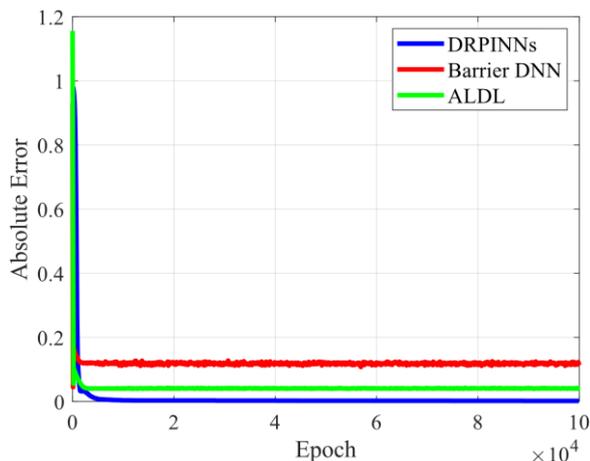

Figure 16 Error comparison for DRPINNs, Barrier DNN, and ALDL on Example 3

In Example 3, the error evolutions of the three methods are shown in Figure 16. The DRPINNs method quickly converges during the early stages of training and stabilizes at an extremely low error level, demonstrating superior accuracy and stability. In contrast, the ALDL method reaches a certain level of stability after an initial decrease but still results in a bigger error than DRPINNs, while the Barrier DNN method yields the largest error. Overall, DRPINNs continues to exhibit the best numerical performance.

In summary, DRPINNs consistently demonstrate superior performance through the three test problems, with faster convergence, lower error, and stronger stability, compared with Barrier DNN method and ALDL method.

## IV  CONCLUSION

We construct a DRPINNS to approximate the variational inequality problems in this paper. First, we transform the variational inequality problem into an optimization problem by using Ritz variational method, then search the optimal loss function weights by Bayesian optimization, and update the training dataset by residual-based adaptive dataset update strategy. Numerical experiments demonstrate that the proposed method outperforms Barrier DNN and ALDL in terms of adaptability and accuracy.